\documentclass[12pt,a4paper,draft]{article}
\usepackage{latexsym,amssymb,amsmath}
\usepackage{ucs}
\usepackage[cp1251,utf8x]{inputenc}
\usepackage[T2A]{fontenc} 
\usepackage[russian,romanian,english]{babel}

\title{On parametrical expressibility in the free void-generated
diagonalizable
algebra\thanks{\copyright Andrei Rusu}}

\author{
Andrei RUSU\\
{\small Information Society Development Institute, andrei.rusu@idsi.md}}

\date{\ }

\begin{document}

\maketitle

Let $\mathfrak A$ be any universal algebra.
We say formula $A$ is explicitly expressible on algebra $\mathfrak A$ via
system of formulas $\Sigma$, if $A$ can be obtained on $\mathfrak A$
of variables and
formulas of $\Sigma$ by means of superpositions. We say a system
of formulas $\Sigma$ is complete in $\mathfrak A$, if any formula is
expressible via $\Sigma$. We say a system $\Sigma$ is precomplete as
to expressibility on $\mathfrak A$ if $\Sigma$ is not complete on
$\mathfrak A$, but for any formula $F$, which is not expressible via
$\Sigma$ on $\mathfrak A$, then the system $\Sigma\cup\{F\}$ is complete
as to expressibility on $\mathfrak A$.
It is known \cite{Post21,Post41} that there are only five precomplete
as to explicite expressibility classes of boolean functions, that there
are only finitely many precomplete classes of functions in
any general $k$-valued logic \cite{Rosenberg65}, that there
are only 12 precomplete classes of pseudo-boolean functions
\cite{Rata71,Rata82}, and other similar results. Note that
preudo-boolean functions cannot be defined by finite tables.

At the same time we can consider other tools to get new functions
of a given system of functions.
We say formula $A$ is parametrical expressible on algebra
$\mathfrak A$ via system of formulas $\Sigma$ if there exist numbers
$l$ and $m$, variables $\pi, \pi_1, \dots, \pi_l$, not occuring in $A$,
formulas $B_1, C_1, \dots, B_m, C_m$, which are explicitly expressible
on $\mathfrak A$ via $\Sigma$, and formulas $D_1 \dots, D_l$ such
that next relations
are valid on $\mathfrak A$:
$$
   (A=\pi) \Longrightarrow
             (\mathop{\wedge}\limits_{i=1}^{m}) (B_i=C_i)
             [\pi_1/D_1] \dots [\pi_l/D_l],
$$
$$
   (\mathop{\wedge}\limits_{i=1}^{m}) (B_i=C_i) \Longrightarrow
   (A=\pi).
$$
In the case of $k$-valued logics this definition was given by
A.V. Kuznetsov as was reported in \cite{Danilcenco81}, and he also
had proved that a two-element set has only 25 parametrically
closed classes of functions, A.F. Danil'\u{c}enko \cite{Danilcenco77}
subsequently proved that a three-element set has only finitely many
parametrically closed classes of functions, S. Burris and R. Willard
\cite{Burris-Willard87} proved recently that a $k$-element set also
has only finitely many parametrically closed classes of functions.
The main result of this paper is that there are systems of functions,
namely of functions of diagonalizable algebra,
such that contains infinitely many parametrically closed classes of
functions.

A diagonalizable algebra $\mathfrak D$ is a boolean algebra
$\mathfrak A = (A; \&, \vee, \supset, \neg)$ with an additional
operator $\Delta$ satisfying the following identies:  $$
  \Delta(\alpha\supset\beta)\le\Delta\alpha\supset\Delta\beta,
$$
$$
  \Delta\alpha\le\Delta\Delta\alpha,
$$
$$
  \Delta(\Delta\alpha\supset\alpha)=\Delta\alpha,
$$
$$
  \Delta 1= 1,
$$
where $1$ is the unit of $\mathfrak A$.

We consider the diagonalizable algebra
$\mathfrak M^* = (M; \&, \vee, \supset, \neg,
\Delta)$ of all infinite binary sequences of the type
$\alpha = (\mu_1, \mu_2, \dots), \mu_i\in\{0, 1\}, i=1,2,\dots$.
The boolean operations $\&, \vee, \supset, \neg$ over elements of $M$
are defined component componentwise, and the operation $\Delta$ over
element $\alpha$ we define by the equality
$\Delta\alpha = (1, \nu_1, \nu_2, \dots)$, where
$\nu_i=\mu_1\&\dots\&\mu_2$. We consider then the subalgebra
generated $\mathfrak M^*$ of the algebra $\mathfrak M$ which is
generated by its zero element $(0, 0, \dots)$.

We prove next

\emph{THEOREM. There are infinitely many precomplete with respect to
parametrical expressibility classes of functions in the free
diagonalizable algebra $\mathfrak M^*$.}

The theorem is based on the example of an infinite family of
parametrically precomplete classes of formulas presented in the
following.

\emph{EXAMPLE. The classes $K_1$, $K_2, \dots$ of formulas, which
preserve on algebra $\mathfrak M^*$ respectively the relations
$x = \neg\Delta 0$, $x = \neg\Delta^2 0, \dots$, constitute a
numerable collection of parametrically precomplete in $\mathfrak M^*$
classes of formulas.}

It is known \cite{Danilcenco81} that these classes are closed
with respect to parametrical expressibility.

It is easy to check that functions $\neg\Delta^i 0$,
$p\& q$, $p\vee q$ belong to the class $K_i$,
and $\neg p$, $\Delta p$ does not belong to $K_i$, $i=1, 2, \dots$.
So these classes are not complete as to parametrical expressibility.

Let us show that they are distimct two by two. It is clear that
the function $\neg\Delta^j 0\not\in K_i$, $i\ne j$.

Let us prove that these classes are parametrically precomplete.
Suppose we have an arbitrary function $F_i(p_1, \dots, p_n)\not\in K_i$.
It means that
$F_i(\neg\Delta^i 0, \dots, \neg\Delta^i 0)\ne \neg\Delta^i 0$.
Let us denote by $\square p$ the function $p\&\Delta p$, and by
$\nabla p$ the function $\square\neg\square\neg\square p$.

Let us consider in the following two functions, denoted respectively
by $F_\neg$ and $F_\Delta$:
$$
  ((\nabla\neg(p\sim q)\&((\neg p\sim q)\sim
    F_i(\neg\Delta^i 0, \dots, \neg\Delta^i 0))) \vee
  (\nabla(p\sim q)\&\neg\Delta^i 0),
$$
$$
  (\nabla q\&((\Delta p\sim q)\sim
    F_i(\neg\Delta^i 0, \dots, \neg\Delta^i 0))) \vee
  (\neg\nabla q\&\neg\Delta^i 0).
$$
It is clear that functions $F_\neg$ and $F_\Delta$ are from class 
$K_i$. 
Then it is not so difficult to check that next relations are valid
on $\mathfrak M^*$:
$$
   (\neg p=q) \Longleftrightarrow
       (F_\neg(p,q)=F_i(\neg\Delta^i 0, \dots, \neg\Delta^i 0)),
$$
$$
   (\Delta p=q) \Longleftrightarrow
       (F_\Delta(p,q)=F_i(\neg\Delta^i 0, \dots, \neg\Delta^i 0)).
$$
The last thing together with previous facts means that classes
$K_1, K_2, \dots$ are parametrically precomplete on
$\mathfrak M^*$.
\par The theorem is proved.

{\sc Remark.} The theorem was announced for the first time at the CAIM-1999 organized by ROMAI (Romanian Society of Applied and Industrial Mathematics) at Pitești, Romania.


\begin{thebibliography}{x}
\bibitem{Post21}
Post E.L. \emph{Introduction to a general theory of elementary
propositions} // Amer. J. Math., 1921, v. 43, p. 163--185.

\bibitem{Post41}
Post E.L. \emph{Two-valued iterative systems of mathematical logic.}
Princeton, 1941.

\bibitem{Rosenberg65}
Rosenberg I. \emph{La structure des functions de plusieure
variables sur un ensemble finit} // C.R. Acad. Sci., 1965, v. 260,
p. 2817--2819.

\bibitem{Rata71}
Ratsa M.F. \emph{The criterion of functional completeness in
the intuitionistic propositional logic} // Dokl. AN SSSR,
1971, vol.  201, no. 4, p. 794--797 (in russian).

\bibitem{Rata82}
Ratsa M.F. \emph{On functional completeness in the intuitionistic
propositional logic} //
Problems of Cybernetics, 1982, vol. 39, p. 107--150 (in russian).

\bibitem{Danilcenco81}
Danil'\u{c}enco A. F. \emph{On parametrical expressibility of
the functions of $k$-valued logic} // Colloq. Math. Soc. Janos
Bolyai, vol. 28, North-Holland, 1981, pp. 147--159.

\bibitem{Danilcenco77}
Danil'\u{c}enco A. F. \emph{Parametric expressibility of
functions of three-valued logic} // Algebra i Logika, vol. 16
(1977), pp.397--416 (in russian)

\bibitem{Burris-Willard87}
Burris S., Willard R. \emph{Finitely many primitive positive clones}
// Proceedings of the American Mathematical Society, vol. 101, no. 3,
1987, pp. 427--430.
\end{thebibliography}
\end{document}